\documentclass[a4paper]{article}
\pdfoutput=1
\usepackage{amsmath}
\usepackage{stmaryrd}
\usepackage{graphicx}
\usepackage{amssymb}
\usepackage{algorithm}
\usepackage{booktabs}
\usepackage{multirow}
\usepackage{algpseudocode}

\usepackage{subfigmat}
\usepackage{bm}
\usepackage{amsmath,amsfonts,mathtools}
\usepackage[width=7in]{geometry}
\providecommand{\mat}[1]{\bm{#1}}%
\renewcommand{\vec}[1]{\mathbf{#1}}

\providecommand{\mA}{\ensuremath{\mat{A}}}

\providecommand{\mD}{\ensuremath{\mat{D}}}

\providecommand{\mI}{\ensuremath{\mat{I}}}

\providecommand{\mQ}{\ensuremath{\mat{Q}}}
\providecommand{\mR}{\ensuremath{\mat{R}}}

\providecommand{\vb}{\ensuremath{\vec{b}}}

\providecommand{\vx}{\ensuremath{\vec{x}}}

\title{
    Kronecker Product Least Squares
}

\author{Pranay Seshadri \thanks{Research Associate, Department of Engineering, University of Cambridge, U.K. }}

\begin{document}

\maketitle

\begin{abstract}
In this rather brief note we present and discuss techniques for solving Kronecker matrix product least squares problems. Our main contribution is an iterative approach that uses the efficient Kronecker matrix-vector multiplication strategy in \cite{fernandes1998efficient} with a conjugate gradient solver. Numerical results contrast this approach---in terms of running times and accuracy---against the direct approach.
\end{abstract}

\section{Introduction}
Kronecker products routinely appear in uncertainty quantification \cite{smith2013uncertainty}, signal processing, image processing and more recently in quantum computing \cite{van2000ubiquitous} and network science \cite{leskovec2010kronecker}. Their algebraic identities enable scalable and rapid matix computations, powering nuermous algorithms and solution techniques. Chief among these is the solution of linear systems, where Kronecker product identities have been applied to solving least squares \cite{fausett1994large}, constrained least squares (using the null space method)\cite{barrlund1998efficient}, stacked least squares (using the generalized singular value decomposition; see page y of \cite{golub2012matrix}) and weighted least squares \cite{coleman2003segmentation}. Our objective in this manuscript is to offer solution strategies for the \emph{multilinear} Kronecker product least squares problem
\begin{equation}
\underset{x}{min}\left\Vert \left(\mA_{1}\otimes \mA_{2}\otimes\ldots\otimes \mA_{n}\right) \vx-\vb\right\Vert _{l_2}, \; \;
\label{eq:main}
\end{equation}
where we assume that $\mA_{k}$ for $k=1, \ldots, n$ are full rank square matrices. It is imperative that during the solution process to this linear system the full Kronecker product is never computed. While this problem has been addressed in literature for the case where $n=2$, to the best of our knowledge no solution exists for the more general case where $n >2$. 

\subsection{Notation and Preliminaries}
We now introduce some notation and few Kronecker product properties that will be useful. Let $\mA_{1} \in \mathbb{R}^{n \times n}$ and $\mA_{2}\in\mathbb{R}^{m\times m}$, then their Kronecker product is
\begin{equation}
\mD = \mA_{1} \otimes \mA_{2},
\end{equation}
where is a $\mD \in \mathbb{R}^{mn \times mn}$ matrix. Essentially each block of $\mD$ is formed by the product of $\mA_{1}$ with the corresponding block of $\mA_{2}$. 

Kronecker products obey algebraic bilinearlity and associativity rules. We enumerate two key properties of Kronecker products that will be useful moving forwards.
\begin{enumerate}
\item The \emph{mixed product property}: For four matrices $\mA_{1}, \mA_{2}, \mA_{3}$ and $\mA_{4}$
\begin{equation}
\left(\mA_{1}\otimes \mA_{2}\right)\left(\mA_{3}\otimes \mA_{4}\right)=\mA_{1}\mA_{3}\otimes \mA_{2}\mA_{4},
\end{equation}
assuming the sizes of the matrices permit the right hand side products to be formed. This property is especially useful in QR factorizations;
\item The \emph{Kronecker matrix vector product}: For a vector $\vx \in \mathbb{R}^{mn}$, the property
\begin{equation}
\left(\mA_{1}\otimes \mA_{2}\right) \vx  = \mA_{2}\mathsf{mat}(\vx) \mA_{1}^{T}
\end{equation}
holds. The operator $\mathsf{mat}(\vx)$ refers to a $\mathbb{R}^{m \times n}$ matrix comprising of all the $nm$ elements in the vector $\vx$. Thus, one can interpret $\vx$ as a vector of the stacked columns of $\mathsf{mat}(\vx)$. 
\end{enumerate}
These two rules are used to develop an expression for solving the Kronecker product least squares problem for the special case where $n=2$ in \eqref{eq:main}.

\subsection{The Direct Least Squares Solution}
Equation \eqref{eq:main} may be re-written in the form of the normal equations
\begin{equation}
\left(\mA_{1}\otimes \mA_{2}\right)^{T}\left(\mA_{1}\otimes \mA_{2}\right) \vx=\left(\mA_{1}\otimes \mA_{2}\right)^{T} \vb.
\end{equation}
Substituting matrices $\mA_1$ and $\mA_2$ with their \emph{thin QR} factorizations gives us
\begin{equation}
\left(\mR_{1}^T\mQ_{1}^T\otimes \mR_{2}^T\mQ_{2}^T\right)\left(\mQ_{1}\mR_{1}\otimes \mQ_{2}\mR_{2}\right)\vx=\left(\mR_{1}^T\mQ_{1}^T\otimes \mR_{2}^T\mQ_{2}^T\right)\vb.
\label{eq:interm}
\end{equation}
Now using the mixed-product property, \eqref{eq:interm} simplifies to
\begin{equation}
\left(\mR_{1}^{T}\mR_{1}\otimes \mR_{2}^{T}\mR_{2}\right) \vx=\left(\mR_{1}^T\mQ_{1}^T\otimes \mR_{2}^T\mQ_{2}^T\right)\vb
\end{equation}
\begin{equation}
\mR_{2}^{T}\mR_{2} \mathsf{mat}(\vx) \mR_{1}^{T} \mR_{1} = \mR_{2}^T \mQ_{2}^T \mathsf{mat}(\vb) \mQ_{1} \mR_{1},
\end{equation}
and since $\mR_{1}$ and $\mR_{2}$ are non-singular, this yields
\begin{equation}
\mR_{2} \mathsf{mat}(\vx)  \mR_{1}^{T}  =  \mQ_{2}^T \mathsf{mat}(\vb)  \mQ_{1},
\end{equation}
\begin{equation}
\Rightarrow \mathsf{mat}(\vx)   =  \mR_{2}^{-1} \mQ_{2}^T \mathsf{mat}(\vb)  \mQ_{1} \mR_{1}^{-T}.
\end{equation}
Similar solution techniques may be achieved using a QR column pivoting factorization and the SVD (singular value decomposition) as discussed in \cite{fausett1994large}. Extending this approach to higher orders remains a challenge.

\section{An Iterative Approach}
In light of the above, we approach this problem using an iterative solution, which can in theory be generalized to any $n$. Our idea is to tie a conjugate gradient solver (see page 628 of \cite{golub2012matrix}) with the effective matrix-vector product multiplication strategy in Fernandes et al.\cite{fernandes1998efficient}; their strategy seeks to compute $\vx \left(\otimes_{k=1}^{n}\mA_{k}\right)$ without ever forming the full tensor product. Their algorithm is based on the tensor product decomposition
\begin{align}
\vx \left(\otimes_{k=1}^{n}\mA_{k}\right) & =  \vx\prod_{k=1}^{n}\mI_{m_{1}}\otimes\ldots\otimes \mI_{m_{k-1}}\otimes \mA_{k}\otimes \mI_{m_{k+1}}\otimes\ldots\otimes \mI_{m_{n}} \\
& = \vx \prod_{k=1}^{n} \mI_{(1:k-1)} \otimes \mA_{k} \otimes \mI_{(k+1:n)},
\end{align}
where $\mI$ are identity matrices. The notation $\mI_{(1:k-1)}$ refers to an identity matrix with $\prod_{g=1}^{k-1}n_{g}$ rows and $\mI_{(k+1:n)}$ to an identity matrix with $\prod_{g=k+1}^{n}n_{g}$ rows.


\section{Numerical Results}
We apply the iterative approach to the problem of approximating the coefficients of tensor grid polynomials. Our codes were implemented in MATLAB and utilized the software's vanilla conjugate gradient solver. The experiments below were carried out on a 8GB Mac with a 2.8 GHz Intel Core i5 processor.

\begin{figure}
\begin{center}
\includegraphics[scale=0.5]{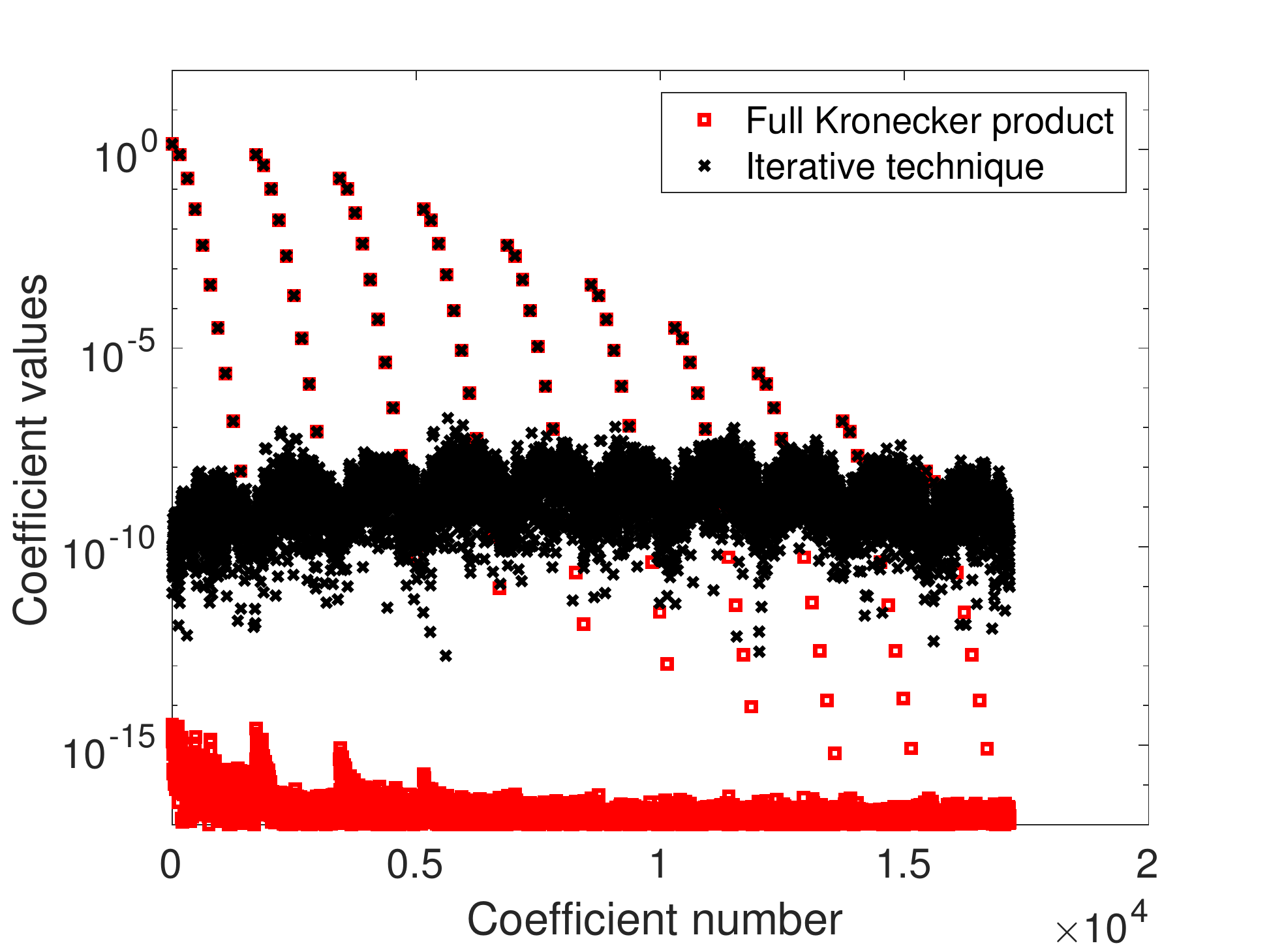}
\caption{Coefficient comparison between the direct method and our iterative approach for the case where $(\mA_1, \mA_2, \mA_3, \mA_4)$ has sizes $(10\times10, 11\times 11, 12\times 12, 13\times 13)$ respectively.}
\end{center}
\label{convresults}
\end{figure}

Table~\ref{tab} presents representative results for the case where $n=4$ for different sizes for each of the matrices (which were all square). The table lists the total runtime for our iterative approach compared to that of the direct approach---which requires computing the full Kronecker product prior to solving the least squares problem. 

In general it was observed that as the matrices get larger, the iterative technique proves to be more efficient from the perspective of run time. To evaluate the numerical accuracy of our approach, Figure~\ref{convresults} plots the least squares coefficients from a direct computation and from the iterative technique for case no. 4 in the table. Coefficient values above $10^{-8}$ (which was our tolerance limit for conjugate gradient solver), were found to match those of the direct solve.

\begin{table}
    \caption{A comparison of the conjugate gradient solver and the full Kronecker product}
\begin{center}
    \begin{tabular}{l l l l l}
    \hline
    No. & Size $(\mA_1, \mA_2, \mA_3, \mA_4)$ & Iterative method runtime (seconds) & Direct method runtime (seconds) \\ \hline
	1 & 2,3,4,5 & 0.33 & 0.004 \\
	2 & 6,7,8,9 & 2.43 & 0.62 \\
	3 & 8,9,10,11 & 5.74 & 6.27 \\
	4 & 10,11,12,13 & 8.14 & 51.47  \\
    \hline
    \end{tabular}
\end{center}
\label{tab}
\end{table}

\end{document}